\documentclass[11pt]{article}%
\usepackage{amsmath}
\usepackage{amsfonts}
\usepackage{amssymb}
\usepackage{graphicx}%
\providecommand{\U}[1]{\protect\rule{.1in}{.1in}}

\begin{document}

\title{On a formula of Daskalopoulos, Hamilton and Sesum}
\author{Bennett Chow}
\date{}
\maketitle

{\scriptsize \textquotedblleft Picture paragraphs unloaded, wise words being
quoted.\textquotedblright\ From `Hail Mary' by Tupac Shakur\smallskip}

We give an exposition of a formula proved by Daskalopoulos, Hamilton and Sesum
\cite{DHS-AncientSurfaceKingRosenau}, which is one of several estimates which
were used in \cite{DHS-AncientSurfaceKingRosenau} to prove that an ancient
solution of the Ricci flow on $\mathcal{S}^{2}$ must be either round or the
King--Rosenau sausage model (\textbf{KR}) solution (see \cite{FOZ},
\cite{King}, \cite{Rosenau} for the KR).

Let $g(t)=\frac{1}{v}g_{\mathcal{S}^{2}}$ be a solution to Ricci flow, where
$g_{\mathcal{S}^{2}}$ is the standard metric on $\mathcal{S}^{2}$. The scalar
curvature is $R_{g}=\frac{1}{v}\frac{\partial v}{\partial t}=\Delta
v-\frac{\left\vert \nabla v\right\vert ^{2}}{v}+2v$, where $\Delta$, $\nabla$
and $\left\vert \cdot\right\vert $ are all with respect to $g_{\mathcal{S}%
^{2}}$. Define $b\doteqdot\operatorname{S}(\nabla^{3}v)$, where
$\operatorname{S}(\alpha)\left(  X,Y,Z\right)  \doteqdot\frac{1}{3}%
(\alpha\left(  X,Y,Z\right)  +\alpha\left(  Y,Z,X\right)  +\alpha\left(
Z,X,Y\right)  )$. The trace-free part is%
\[
\operatorname{TF}(b)\left(  X,Y,Z\right)  \doteqdot b\left(  X,Y,Z\right)
-\left(  z\left(  X\right)  \left\langle Y,Z\right\rangle +z\left(  Y\right)
\left\langle Z,X\right\rangle +z\left(  Z\right)  \left\langle
X,Y\right\rangle \right)  ,
\]
where $z=\frac{1}{4}\operatorname{tr}_{g_{\mathcal{S}^{2}}}^{2,3}b=\frac{1}%
{4}(d\Delta v+\frac{2}{3}dv)$. Given a $4$-tensor $c$ symmetric and trace-free
in the last $3$ slots,%
\begin{align*}
c\left(  W,X,Y,Z\right)   &  \doteqdot\overset{\circ}{c}\left(
W,X,Y,Z\right)  +\left\langle W,X\right\rangle e\left(  Y,Z\right)
+\left\langle W,Y\right\rangle e\left(  Z,X\right)  +\left\langle
W,Z\right\rangle e\left(  X,Y\right) \\
&  \quad\;+\left\langle Y,Z\right\rangle f\left(  X,W\right)  +\left\langle
Z,X\right\rangle f\left(  Y,W\right)  +\left\langle X,Y\right\rangle f\left(
Z,W\right)  ,
\end{align*}
where $e=\frac{1}{3}\operatorname{tr}_{g_{\mathcal{S}^{2}}}^{1,2}c=-2f$
($\overset{\circ}{c},e,f$ are totally trace free). Then $Q\doteqdot
v\left\vert \operatorname{TF}(b)\right\vert ^{2}$ satisfies ($L\bar{Q}$ in
\cite[\S 5]{DHS-AncientSurfaceKingRosenau})%
\[
\frac{\partial}{\partial t}Q=v\Delta Q-4RQ-2\left\vert \overset{\circ
}{v\nabla\operatorname{TF}(b)+2dv\otimes\operatorname{TF}(b)}\right\vert
^{2}-\frac{1}{2}\left\vert v\operatorname{TF}\left(  \nabla^{2}(\Delta
v+6v)\right)  -2\operatorname{tr}{}_{g}^{1,2}\left(  dv\otimes
\operatorname{TF}(b)\right)  \right\vert ^{2},
\]
where $\operatorname{TF}\left(  \alpha\right)  \doteqdot\alpha-\frac{1}%
{2}(\operatorname{tr}_{g_{\mathcal{S}^{2}}}\alpha)g_{\mathcal{S}^{2}}$ for a
symmetric $2$-tensor $\alpha$. Note that $Q$ vanishes on the KR
solution.\smallskip

\textbf{Remark.} Taking $\varphi=\log v$ in $\Delta_{g}(R+\left\vert
\nabla\varphi\right\vert _{g}^{2})=A+2g\left(  \nabla(\Delta_{g}%
\varphi-R),\nabla\varphi\right)  +\left(  \Delta_{g}\varphi\right)  ^{2}%
-R^{2}$, where $A\doteqdot\Delta_{g}R+R^{2}-\frac{\left\vert \nabla
R\right\vert _{g}^{2}}{R}+\frac{\left\vert \nabla R+R\nabla\varphi\right\vert
_{g}^{2}}{R}+2|\nabla_{g}^{2}\varphi-\frac{1}{2}\Delta_{g}\varphi g|_{g}%
^{2}\geq0$ by the trace Harnack estimate for ancient solutions, we get
$\Delta_{\mathcal{S}^{2}}\left(  \Delta_{\mathcal{S}^{2}}v+6v\right)
\geq-4v\geq-C$ (\cite[(2.7)]{DHS-AncientSurfaceKingRosenau}). Define
$J_{\alpha}\left(  g\right)  =\int_{\mathcal{S}^{2}}(\frac{\left\vert \nabla
v\right\vert _{\mathcal{S}^{2}}^{2}}{v^{\alpha}}+F_{\alpha}\!\left(  v\right)
)d\mu_{\mathcal{S}^{2}}$, where $F_{\alpha}\!\left(  v\right)  =-\frac
{4}{2-\alpha}v^{2-\alpha}$ ($\alpha\neq2$), $F_{2}\!\left(  v\right)  =-4\log
v$. Then $\frac{d}{dt}J_{\alpha}\left(  g\left(  t\right)  \right)
=\int_{\mathcal{S}^{2}}\left(  -2\frac{\partial v}{\partial t}-\left(
2-\alpha\right)  \left\vert \nabla v\right\vert _{\mathcal{S}^{2}}^{2}\right)
\frac{\frac{\partial v}{\partial t}}{v^{\alpha+1}}d\mu_{\mathcal{S}^{2}}$,
which is $\leq0$ if $\alpha\leq2$. $J_{2}$ is Polyakov's energy; $J_{1}$ is
\cite[(3.4)]{DHS-AncientSurfaceKingRosenau}.

\end{document}